\documentclass[final,epsfig,amsfont]{article}
\usepackage{latexsym}
\usepackage[dvips]{graphicx}
\newcommand{\eeq}{\end{equation}}
\newcommand{\beq}{\begin{equation}}

\begin{document}
\title{Computing the equilibrium measure of a system of intervals converging to a Cantor set}
\author{
Giorgio Mantica \\
Center for Non-linear and Complex Systems,\\
Dipartimento di Scienze ed Alta Tecnologia, \\ Universit\`a dell' Insubria,\\
via Valleggio 11, 22100
Como, Italy. \\ Also at I.N.F.N. sezione di Milano and CNISM unit\`a di Como.
}
\date{}
\maketitle
\begin{abstract}
We describe a numerical technique to compute the equilibrium measure, in logarithmic potential theory, living on the attractor of Iterated Function Systems composed of one-dimensional affine maps. This measure is obtained as the limit of a sequence of equilibrium measures on finite unions of intervals. Although these latter are known analytically, their computation requires the evaluation of a number of integrals and the solution of a non-linear set of equations. We unveil the potential numerical dangers hiding in these problems and we propose detailed solutions to all of them. Convergence of the procedure is illustrated in specific examples and is gauged by computing the electrostatic potential.
\end{abstract}



\begin{center}
{\large \em In memory of Professor Don Luigi Verga} \\
\end{center}

\section{Introduction}
\label{sec1}

Suppose that a finite amount of ``charge'' is placed on a set of ``conductors'' of arbitrary shape, placed at equally arbitrary positions in space, and connected by thin, conducting wires. Then, this charge distributes itself so to minimize repulsion, {\em i.e.}, the electrostatic energy of its configuration. It is very easy for us to intuitively perceive this phenomenon, and indeed we have been exposed to it since early in school, via easily realizable laboratory experiments that go back to the years of Alessandro Volta and to his electrophorus \cite{volta}, which he did not invent, but he perfected into an instrument to accumulate and transfer electric charge, years before his monumental Pila. This intuition is common knowledge among scientists of all disciplines: in Italy, at least until recently, electrostatics was taught in  high school, both scientific and of classical kind and, even more amazingly, it was included in the syllabus of the degree in mathematics and of medical school alike! This paper is dedicated to a wonderful teacher who let middle--school students, twelve years--old, play with electrostatic devices and inspired them to pursue scientific research.

What is far less known is the fact that the mathematical scheme behind this physical description can be easily generalized to spaces and conductors of arbitrary dimension---even fractal \cite{str0,str0b}---and to arbitrary analytical forms of the ``electrostatic'' interaction, that need not be ``Coulomb'', but can take on different analytical forms.
These generalizations have proven to be extremely fruitful in a wide variety of fields, from PDE's to harmonic analysis, to name just a few.
In this paper, we will put ourselves in the situation where the ambient space is two--dimensional (the complex plane) and the pseudo--Coulomb law depends on the inverse of the distance between two charges: physicists will easily think of the interaction of two infinite, parallel wires, mathematicians will immediately recognize the logarithmic potential that takes us into the realm of the well established {\em logarithmic potential theory} \cite{ran0,ed}.

Indeed, if we let $E$ to be a compact domain in the complex plane $\bf C$, and $\sigma$ a positive probability measure supported on $E$, the potential $V(\sigma;z)$ generated by $\sigma$ at the point $z$ in $\bf C$, is given by the formula
 \begin{equation}
 \label{pote1}
    V(\sigma;z) := - \int_E \log |z-s| \; d \sigma(s).
 \end{equation}
The electrostatic energy ${\cal E}(\sigma)$ of the distribution $\sigma$ is then given by the integral
 \begin{equation}
 \label{pote2}
  {\cal E}(\sigma) := \int_E  V(\sigma;u) \; d\sigma(u)  = - \int_E \int_E \log |u-s| \; d\sigma(s)d\sigma(u).
 \end{equation}
The equilibrium measure $\sigma_E$ associated with the compact domain $E$ is the unique measure that minimizes the energy ${\cal E}(\sigma)$, when this latter is not identically infinite. Potential theory gives us a formidable set of technical instruments \cite{ran0,ed} to deal with this problem in relation to general compact sets $E$.

In this paper, we are interested in a specific family of compact sets $E$, that are Cantor sets on the real line associated with the construction of Iterated Function Systems (IFS)---to be defined momentarily---and we study the efficient numerical computation of their equilibrium measure $\sigma_E$ and of its potential $V(\sigma;z)$. Our results will consist of a reliable algorithms to achieve these goals. This algorithm is rather straightforward, in the sense that it combines well known theoretical facts (the analytical form of the equilibrium measure when $E$ is a collection of intervals) and well known numerical techniques (Gaussian summation, root finding) with a standard idea in IFS construction: the hierarchical approach. Despite or rather because of this simplicity, numerical experiments indicate that the algorithm is stable and it can be used as a microscope to probe deeply into the structure of the equilibrium measure on Cantor sets.

This paper is organized as follows: in the next section we define the attractors of Iterated Function Systems, via a hierarchical construction. Each step in this construction yields a finite family of intervals, whose equilibrium measure is analytically known. We sketch this solution in Sect. \ref{sec-harm}.
The following three sections contain the numerical techniques to find this solution and to overcome a few potentially destructive numerical difficulties. In Sect. \ref{sec-integ} we apply  Gaussian integration to compute a class of integrals that reduce the problem to a set of non--linear equations, and in Sect. \ref{sec-solve} we describe a convenient approach for the solution of these equations.
Next, in Sect. \ref{sec-hierar}, we show how to tailor the root finding routine to the hierarchical structure of IFS attractors. As a result, we obtain the equilibrium measure on these latter. Pictures and graphs of numerical experiments illustrate the theory. In section \ref{sec-capa} we compute the potential generated in the complex plane by the equilibrium measure on the IFS attractor and we compute its capacity.  The application of this approach to deep questions on equilibrium measures on Cantor sets is briefly discussed in the conclusions.

\section{Attractors of Iterated Function Systems}
\label{sec-ifs}

Let us therefore construct the set $E$ of which we want to place an electrostatic charge. To do this, we need to recall, as briefly as possible and in the simplest setting, the construction of Iterated Function Systems (IFS) \cite{papmor,hut,diaco,dem,ba2}.
These are collections of maps $\phi_i : {\bf R} \rightarrow
{\bf R}$, $i = 1, \ldots, M$, for which there exists a
set ${\mathcal A}$, called the {\em attractor} of the IFS, that solves
the equation
 \begin{equation}
 \label{attra}
    {\mathcal A}=\bigcup_{j=1,\ldots ,M}\;\phi_i({\mathcal A}) :=  \Phi ({\cal A}).
 \end{equation}
When the maps $\phi_i$ are contractive, the attractor
${\mathcal A}$ is unique. It can also be seen as the fixed point of the operator $\Phi$, defined in eq. (\ref{attra}), on the set of compact subsets of ${\bf R}$. Since this space is complete in the Hausdorff metric, and since $\Phi$ turns out to be contractive in this metric,
the set ${\cal A}$ can be also found as the limit of the sequence $\Phi^{n}(E^{0})$,
where $E^{0}$ is any non-empty compact set:
\begin{equation}
{\cal A} = \lim_{n \to \infty} \Phi^{n}(E^{0}).
\label{eq-itera}
\end{equation}

Attractors of Iterated Function Systems feature a rich variety of topological structures, so that their full characterization is far from being fully understood, especially in the case of IFS with uncountably many maps \cite{nalgo1,nalgo2,arxiv,intjns}. We nonetheless restrict ourselves in this paper to a particular family, that of IFS with a finite number of affine maps of the form:
\begin{equation}
\label{mappi}
    \phi_j (s) = \delta_j (s - \gamma_j) + \gamma_j,  \;\;  j = 1, \ldots, M ,
\end{equation}
where $\delta_i$ are real numbers between zero and one, called {\em contraction ratios}, and $\gamma_{i}$ are real constants, that geometrically correspond to the fixed points of the maps.
Under these conditions, the attractor ${\cal A}$ is a finite or infinite collection of intervals, or a Cantor set. We will consider this last, interesting case. It can be easily recognized that the famous ``middle--third Cantor set'' follows in this class: it can be obtained by taking just two maps with $\delta=1/3$ and $\gamma_1=-1$, $\gamma_2=1$. We shall study this case later in this paper.

As remarked above, any non-empty compact set $E^0$ can be used in eq. (\ref{eq-itera}). We let $E^0$ be the convex hull of the attractor ${\cal A}$, that can be easily be identified as the interval
\begin{equation}
\label{eq-int0}
    E^0 =  [\gamma_1,\gamma_M],
\end{equation}
where we have ordered the IFS maps according to increasing values of their fixed points:  $\gamma_j < \gamma_{j+1}$, for any $j=1,\ldots,M-1$.
This is the core of the common hierarchical construction of the set ${\cal A}$, as the limit of the sequence of compact sets $E^n$:
\begin{equation}
\label{eq-int1}
    E^n = \Phi^n(E^0) = \bigcup_{i=1}^{M^n} [\alpha_i,\beta_i].
\end{equation}
In the case of {\em fully disconnected} IFS ({\em i.e.} those for which the intervals $\phi_j(E^0)$ are pairwise disjoint), $E^n$ is the union of $N:=M^n$ disjoint intervals, that we will denote as $E^n_i := [\alpha_i,\beta_i]$.  In other words, the set $E^n$ can be seen as a vector of intervals, stored via the two vectors of their extreme points. Of course, these latter should also carry a superscript $n$, referring to the generation. Not to overburden the notation, we will leave this superscript implicit, when confusion is not possible. We restrict ourselves to the disconnected IFS case, for the remainder of the paper.

For similarity with the spectral analysis of periodic solids, we will call the intervals $E^n_i$ at r.h.s. of eq. (\ref{eq-int1}) {\em bands} at generation (or level) $n$. It appears from eq. (\ref{eq-int1}) that all bands $E_i^n$ at level $n$ can be obtained by applying the transformations $\phi_j$, $j=1,\ldots,M$, to the bands at generation $n-1$. As a consequence, the length of these intervals, as well as their sum, tend to zero geometrically and the attractor is a Cantor set of null Lebesgue measure.

While none of the bands $E_i^{n-1}$ are also bands at level $n$, the contrary happens for the so--called {\em gaps} between the bands. For these, we will use the notation:
\begin{equation}
\label{eq-gap1}
   G^n := \bigcup_i G^n_i = \bigcup_{i=1}^{M^n-1} (\beta_i,\alpha_{i+1}).
\end{equation}

At level zero we have no bands, while at level one these are in the number of $M-1$. It is then noticeable that {\em all} gaps at level $n-1$ remain gaps at level $n$, and at the same time a set of new gaps is generated from the former: let $H^n$ be the set of ``new'' gaps created at level $n$, so that
\begin{equation}
\label{eq-int2}
 G^n = G^{n-1} \bigcup H^n.
\end{equation}
The latter set, $H^n$, is iteratively constructed as
\begin{equation}
\label{eq-gap2}
  H^{n} = \bigcup_{j=1}^{M}  \phi_j(H^{n-1}),
\end{equation}
starting from $H^1 = (\beta^1_1,\alpha^1_2) \bigcup \ldots \bigcup (\beta^1_{M-1},\alpha^1_M)$.
This seemingly complicated algorithm is nothing more than the straightforward generalization of the construction of the ternary Cantor set by deleting the ``middle third''. This property of the set of gaps is important in our computation of the equilibrium measure on the Cantor set $\cal A$.

\section{Equilibrium measure on a set of intervals}
\label{sec-harm}

Let us now embed $E^n$ in the complex plane, and study its equilibrium measure. The solution of the equilibrium problem for a finite union of $N$ intervals $[\alpha_i,\beta_i]$, $i = 1,\ldots,N$, is well known \cite{widom,sasha1,nuttall,franz0,franz}. Define the polynomial $Y(z)$,
 \begin{equation}
 \label{meas2}
    Y(z) = \prod_{i=1}^{N} (z-\alpha_i)(z-\beta_i),
 \end{equation}
and its square root, $\sqrt{Y(z)}$, as the one which takes real values for $z$ real and large (larger than $\gamma_M$ indeed). Also, let the real number $\zeta_i$ belong to the open interval $(\beta_i,\alpha_{i+1})$, for $i=1,\ldots,N-1$ ({\em i.e.} to the gap $G^n_i$). Define $Z(z)$ as the monic polynomial of degree $N-1$ with roots at all $\zeta_i$'s:
 \begin{equation}
 \label{meas3}
    Z(z) =  \prod_{i=1}^{N-1} (z-\zeta_i).
 \end{equation}
With these premises, there exists a unique set of values $\{\zeta_i, i = 1,\ldots,N-1\}$ that solve the set of coupled, non--linear equations
 \begin{equation}
 \label{meas5}
    \int_{b_i}^{a_{i+1}} \frac{Z(s)}{\sqrt{|Y(s)|}} \; ds = 0, \;\; i = 1,\ldots,N-1.
 \end{equation}
Our first task in the following will be to evaluate numerically these integrals for any given set of values $\{\zeta_i, i = 1,\ldots,N-1\}$. The second task will be to find the solution of this set of equations. This is of paramount importance, for it permits us to find the equilibrium measure of the set $E^n$: for simplicity of notation, we denote it by $\sigma^n := \sigma_{E^n}$:
 \begin{equation}
 \label{meas6}
     d \sigma^n (s) = \frac{1}{\pi} \sum_{i=1}^N \chi_{[\alpha_i,\beta_i]} (s)\;  \frac{|Z(s)|}{\sqrt{|Y(s)|}} \; ds.
 \end{equation}
It is apparent from the previous equation that the measure $\sigma^n$ is absolutely continuous with respect to the Lebesque measure on $E^n$.
Of particular relevance are also the integrals of $\sigma^n$ over the intervals composing $E^n$, that we call the {\em harmonic frequencies}:
 \begin{equation}
 \label{meas8}
     \omega^n_{i} := \frac{1}{\pi} \int_{\alpha_i}^{\beta_i}  \frac{|Z(s)|}{\sqrt{|Y(s)|}} ds.
 \end{equation}
In fact, when these frequencies are rational numbers, of the kind $p_i/N$, with $p_i$ integer, there exists a strict-T polynomial on $E^n$ \cite{franz}, that is, a polynomial with oscillation properties mimicking, and extending, those of the classical Chebyshev polynomials on a single interval.
Also, we will consider the integrated measures $\Omega^n_{m}$:
\begin{equation}
 \label{meas8b}
     \Omega^n_{m} := \sum_{i=1}^m \omega^n_{i},
 \end{equation}
that are the integral of $\sigma^n$ on the interval $[\gamma_1,\beta_m]$.

Finally, the Green function with pole at infinity can be computed as the complex integral
\begin{equation}
 \label{meas9}
    G(z) = \int_{\alpha_1}^z  \frac{Z(s)}{\sqrt{Y(s)}} ds,
 \end{equation}
in which no absolute value appears. The Green function can be used to compute the logarithmic capacity of the set $E^n$, $C(E^n)$, via a real integral
\begin{equation}
 \label{meas10}
    \log(C(E^n)) = \int_{-\infty}^{\alpha_1}  [\frac{Z(s)}{\sqrt{Y(s)}} - \frac{1}{s-(\alpha_1+1)}] ds.
 \end{equation}
Although our technique could permit to evaluate the above integral, we will compute the capacity $C(E^n)$ following a different approach.

\section{Integrating the ratio $Z/\sqrt{Y}$}
\label{sec-integ}
The first step in the determination of the equilibrium measure is the numerical determination of integrals of the kind (\ref{meas5}), or of the kind (\ref{meas8}), that are quite similar. A straightforward technique to achieve this goal is described in this section.
The first integrals are:
 \begin{equation}
 \label{meas8c}
   {\cal K}_i := \frac{1}{\pi} \int_{G_i}  \frac{Z(s)}{\sqrt{|Y(s)|}} \; ds,
 \end{equation}
where $G_i$ is a gap at generation $n$. The case when $G_i$ is replaced by a band $E_i$ can be handled by minor variations of the technique that we are going to describe. For simplicity of notation, since $n$ is fixed, we will omit its mention in the following derivation. A full set of integrals ${\cal K}_i$, for variable $i$, need to be evaluated. Let us now keep the index $i$ fixed and show how to compute any single one of them.

Firstly, consider the change of variables induced by the affine transformation $\psi_i:= x \rightarrow A_i x + B_i$, that maps the interval $G_i$ into $E^0=[\gamma_1,\gamma_M]$.  Without any loss of generality we can assume that $E^0=[-1,1]$.
The coefficients $A_i$ and $B_i$ can be computed by careful bookkeeping from the map parameters $\delta_j$ and $\beta_j$. Even if we feel no need to report the explicit formulae here, they are essential in the numerical implementation that we will describe in the following.
In so doing, the bands $E_m$ and the gaps $G_m$, for $m$ different than $i$, are mapped, via $\psi_i$, to intervals outside $[-1,1]$. Equivalently, the variables $\zeta_m$ are mapped via the same map to new values $\psi_i(\zeta_m)$. New functions $\bar{Z}_i$ and $\bar{Y}_i$ 
are then written as in eqs. (\ref{meas2}) and (\ref{meas3}), now in terms of the transformed $\psi_i(\zeta_k)$, $\psi_i(\alpha_m)$ and $\psi_i(\beta_m)$. In conclusion, the integral ${\cal K}_i$ can be written as
 \begin{equation}
 \label{meas8d}
   {\cal K}_i = \frac{1}{\pi} \int_{-1}^{1}  \frac{\bar{Z}_i(s)}{\sqrt{|\bar{Y}_i(s)|}} ds.
 \end{equation}

Secondly, remark that the product $1-z^2$ appears in the rescaled function $\bar{Y}_i(z)$. Therefore, we can part the full product $\bar{Y}_i(z)$ in two factors:
 \begin{equation}
 \label{meas2b}
    \bar{Y}_i(z) = (1-z^2) \tilde{Y}_i(z),
 \end{equation}
where $\tilde{Y}_i(z)$ is implicitly defined. With these notations, it appears clearly that the integration of $\frac{\bar{Z}_i(s)}{\sqrt{|\bar{Y}_i(s)|}}$ with respect to the Lebesgue measure on $[-1,1]$ can be seen as an integration with respect to the Chebyshev measure ${ds}/{\pi \sqrt{1-s^2}}$. Then, the integral ${\cal K}_i$ can be approximated by a Gaussian summation,
\begin{equation}
 \label{meas5b}
   {\cal K}_i = \int_{-1}^{1} \frac{\bar{Z}_i(s)}{\sqrt{|\tilde{Y}_i(s)|}} \;
    \frac{ds}{\pi \sqrt{1-s^2}} \simeq \sum_{k=1}^{K} w_k \frac{\bar{Z}_i(x_k)}{\sqrt{|\tilde{Y}_i(x_k)|}},
 \end{equation}
where $x_k,w_k$ are Gaussian points and weights for the Chebyshev measure, respectively. As a consequence, each integral ${\cal K}_i$ can be readily evaluated by Gaussian summation.

Clearly, the above assumes that the values $\zeta_i$ be known for all $i=1,\ldots,N-1$.
To compute these latter, the $N$ equations ${\cal K}_i=0$, $i=1,\ldots,N-1$, must be solved. In the next section we employ an algorithm for the solution of this set of non-linear equation that requires the computation of the derivatives $\partial {\cal K}_i / \partial \zeta_m$. Also these quantities can be computed via Gaussian summation. Observe in fact that the partial derivative of eq. (\ref{meas5b}) requires the computation of
 \begin{equation}
 \label{meas12}
    \frac{\partial}{\partial \zeta_m} \bar{Z}_i (z) = - A_i \prod'_{l \neq m} (z-\psi_i(\zeta_l)),
 \end{equation}
that therefore leads to
\begin{equation}
 \label{meas13}
 \frac{\partial {\cal K}_i}{\partial \zeta_m}
     \simeq - A_i \sum_{k=1}^{K} w_k \frac{
     \prod'_{l \neq m} (x_k-\psi_i(\zeta_l))
     }{\sqrt{|\tilde{Y}(x_k)|}}.
 \end{equation}

\section{Solving the non-linear equations ${\cal K}_i =0$}
\label{sec-solve}

We have seen that the determination of the equilibrium measure on the set $E^n$ requires the solution of a system of $N-1$ non--linear functions in $N-1$ variables, ${\cal K}_i (\zeta_1,\ldots,\zeta_{N-1}) =0$.
We can compute this solution by a careful usage of the routine {\sc HYBRJ} in Minpack \cite{minpack}, which employs a modification of the Powell hybrid method. This technique is a variation of Newton's method, that iteratively converges to the solution vector. As it turns out, {\sc HYBRJ} cannot be applied blindly without taking into account the particular nature of the set $E^n$ and of the discretized integrals (\ref{meas5b}).
In fact, in the construction of the attractor of a fully disconnected IFS, all bands $E^n_i$ have a length that goes to zero geometrically fast, while gaps differ in length by orders of magnitude. This requires a series of procedures, in the evaluation of the integrals ${\cal K}_i$ and $\frac{\partial {\cal K}_i}{\partial \zeta_m}$, that we explain in this section.

The first consists in the approach described at the end of the preceding section: it is convenient to scale the each interval $G^n_i$ to $[-1,1]$, when performing the relative integration, to employ the Gaussian form in eqs. (\ref{meas5b}) and (\ref{meas13}).

When this is done, a further difficulty is to be overcome: the function ${\bar{Z}_i(x_k)}/{\sqrt{|\tilde{Y}_i(x_k)|}}$ is made of factors differing by various orders of magnitude, that must be conveniently rearranged when evaluating it numerically. Therefore, in the expansion
\begin{equation}
 \label{meas15}
   \frac{\bar{Z}_i(x_k)}{\sqrt{|\tilde{Y}_i(x_k)|}} =
   \frac{
   \prod_{l=1}^{N-1} (x_k - \psi_i(\zeta_l))}{ \sqrt{
    |x_k-\psi_i(\alpha_i)||x_k-\psi_i(\beta_{m+1})|
    \prod_{m\neq i,i+1} |x_k - \psi_i(\alpha_m)|
   |x_k - \psi_i(\beta_m)|  } }
 \end{equation}
we find it convenient to group together the ratios $(x_k-\psi_i(\zeta_m))/\sqrt{|x_k - \psi_i(\alpha_m)||x_k - \psi_i(\beta_m)|}$, for $m < i$, and
$(x_k-\psi_i(\zeta_m))/\sqrt{|x_k - \psi_i(\alpha_{m+1})||x_k - \psi_i(\beta_{m+1})|}$, for $m > i$, in which numerator and denominator have comparable size. The more the index $m$ is different from $i$, the larger are the differences between $x_k$ and the points $\psi_i(\alpha_m)$ and $\psi_i(\beta_m)$, and the more these ratios tend to one. Physically, this means that remote gaps $G^n_m$ have a small influence on the integral ${\cal K}_i$: our ordering correctly reproduces this fact.
As it appears from eq. (\ref{meas15}), grouping terms as above leaves out a few factors, those close to the interval $\psi_i(G^n_i)=[-1,1]$, that mostly contribute to the value of the full product, and are therefore taken care separately.
A similar procedure can also be applied to the formula for derivatives, eq. (\ref{meas13}).

But the more important precaution to be taken follows from considering the variables involved in the Newton-like iteration of {\sc HYBRJ}. In the original equations ${\cal K}_i =0$ the unknowns are the positions of the $N-1$ zeros of $Z$, {\em i.e.} $\zeta_m$, $m=1,\ldots,N-1$. Each of these belongs to the gap $G^n_m$. It is then apparent, from the description in Sect. \ref{sec-ifs}, that the range of variation of these values can differ by orders of magnitude, and can quickly become smaller than numerical precision. It is then mandatory to resort to new variables: we consider in place of $\zeta_m$ the normalized variables $\lambda_m := \psi_m (\zeta_m)$, where, as above, $\psi_m$ is the affine transformation that maps $G^n_m$ into $[-1,1]$. These new variables are {\em all} bound to the interval $[-1,1]$. Some care must now be taken in the evaluation of the compounded derivatives, but fortunately, the new variables are the most natural, in the construction of the IFS, as it can be easily recognized from eqs. (\ref{mappi}), (\ref{eq-int1}).

\begin{figure}
\centerline{\includegraphics[width=.6\textwidth, angle = -90]{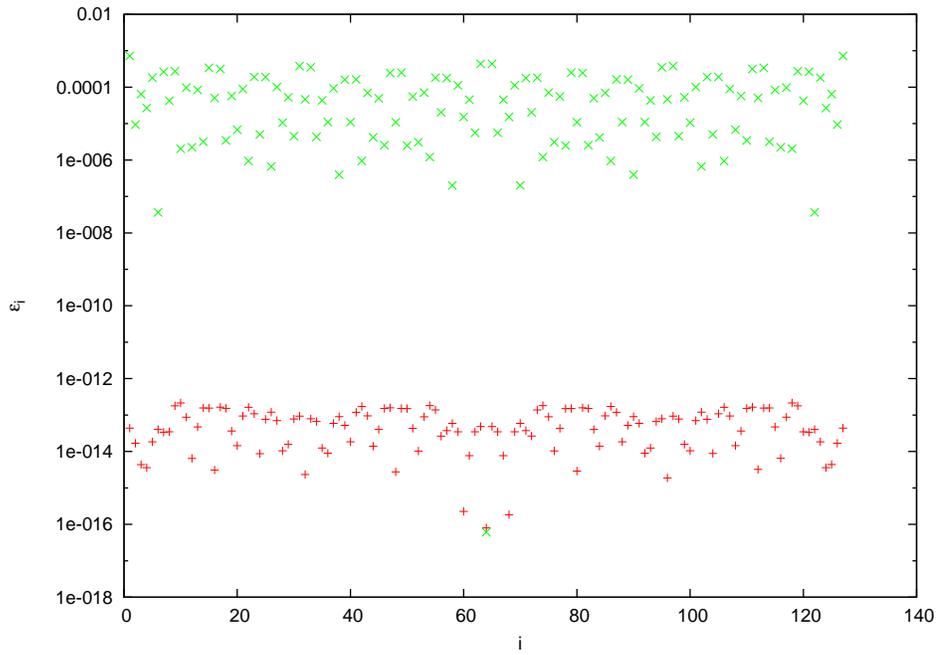}}
\caption{Absolute value $\epsilon_i$ of the non-linear functions ${\cal K}_i$ at the initial vector $(\lambda_1,\ldots,\lambda_{N-1})$ (see next section for its definition) (green crosses), and after the action of the root finding routine (red crosses). The case under study is the ternary Cantor set at generation number $n=7$. The number of Gaussian points employed is 2048.}
\label{erra1}
\end{figure}

The above provides us with a working code for the computation of the non-linear functions (\ref{meas5b}) and of the partial derivatives (\ref{meas12}), as functions of the variables $\lambda_m$, $m=1,\ldots,N-1$, that can be linked as a subroutine to the fortran (yes, fortran!) code {\sc HYBRJ}, run in double precision on a 32-bit processor. In Figure (\ref{erra1}) we show the absolute values of the integrals ${\cal K}_i$ before and after the call to {\sc HYBRJ} in a typical case, that of the ternary Cantor set at generation number $n=7$, to be defined more precisely later on in Sect. \ref{sec-capa}. The algorithm has converged to the correct solution within numerical precision.

\begin{figure}
\centerline{\includegraphics[width=.6\textwidth, angle = -90]{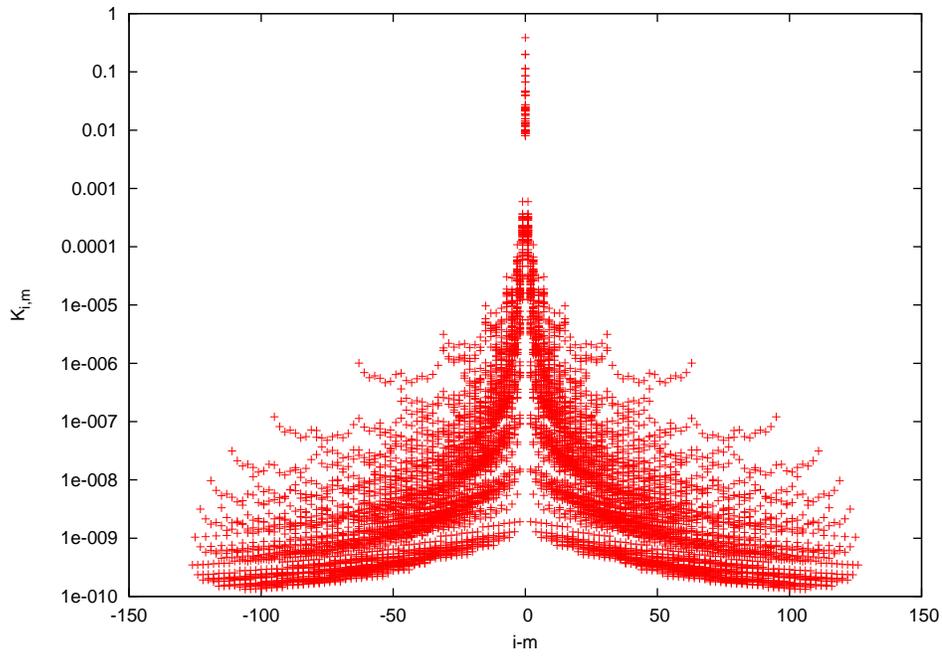}}
\caption{Absolute value $K_{i,m}$ of the partial derivatives $\frac{\partial {\cal K}_i}{\partial \lambda_m}$ versus $i-m$ at the final solution vector, for the same case of Fig. \ref{erra1}.}
\label{erra2}
\end{figure}

To fully appreciate the precision of the solution vector, we can examine the value of the partial derivatives $\frac{\partial {\cal K}_i}{\partial \lambda_m}$ at the solution vector. These values are reported in Figure (\ref{erra2}) versus the difference $i-m$: as remarked above, the larger this difference (in absolute value) the lesser the influence of $\lambda_m$ on the integral ${\cal K}_i$: this is clearly evident in the figure. Also, we can observe that the system of equations is almost ``diagonal'', in the sense that diagonal components with $i = m$ (the tip of the ``Christmas tree'') are orders of magnitude larger than the non-diagonal branches, $i \neq m$.

\section{Hierarchical structure of the equilibrium measure}
\label{sec-hierar}

We have described in the previous section an algorithm for the determination of the roots of the equilibrium equations (\ref{meas5}). The algorithm applies to any finite family of intervals. In particular, we want to apply it to a system of intervals generated by a fully disconnected IFS. In the previous section we have already presented results in Figures \ref{erra1} and \ref{erra2} referring to the case of the middle--third Cantor set. A further remark on the way they have been obtained is important.
While the case of the ternary Cantor set is particularly favorable because of its symmetry, we now choose to describe our procedure in the case of an asymmetric IFS given by two-maps with $\delta_1=4/5$ and $\delta_2=1/10$. We will return to the middle--third Cantor set in the next section.

Observe that the cardinality of $E^n$ is $M^n$, where $M$ is the number of IFS maps, and therefore the number of gaps is $M^n-1$. Because of eq. (\ref{eq-int2}), a subset of cardinality $M^{n-1}-1$ of these latter exist already as gaps at generation $n-1$. The old gaps are approximately $(M-1)$ times less numerous than the new ones, nonetheless, their r\^ole is crucial. In fact, the Newton-like search of the algorithm {\sc HYBRJ} in Minpack need to be initialized. It is particularly efficient to initialize the values of the roots $\lambda^n_{Mm}$ to the ``old''  values at generation $\lambda^{n-1}_m$, and the new ones to zero (which means that $\zeta_m$ lies in the middle of $G^n_m$.)

\begin{figure}
\centerline{\includegraphics[width=.6\textwidth, angle = -90]{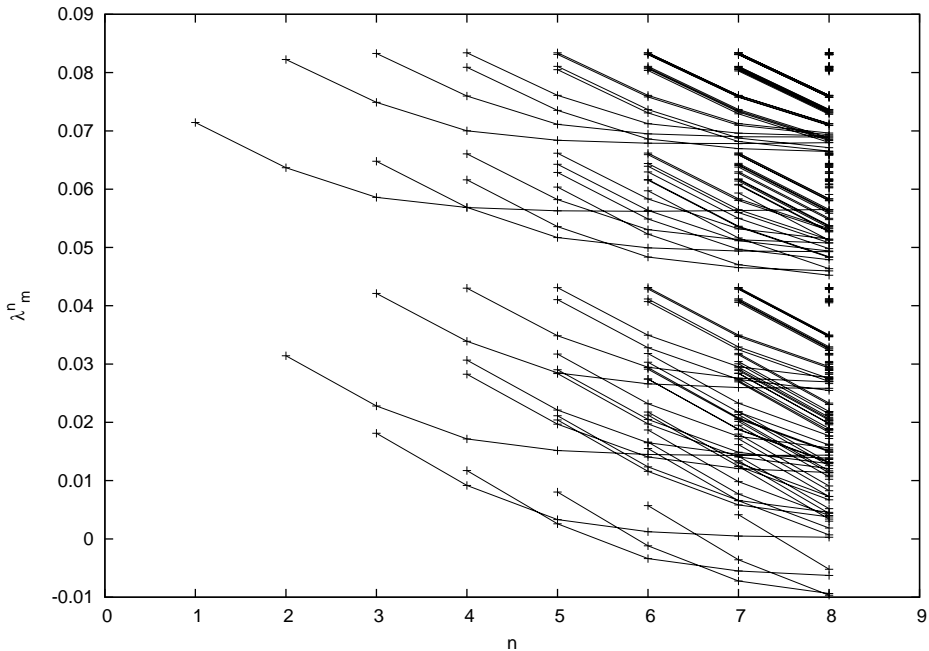}}
\caption{Roots $\lambda^n_m$ of the equations (\ref{meas5}) versus generation number $n$, for $m=1,\ldots,2^{n-1}-1$, for the IFS described in the text. Lines connect roots in the same gap at different $n$. The number of Gaussian points employed is 160.}
\label{porca1}
\end{figure}

Everything described is clearly evident in Fig. \ref{porca1}, that plots the values of the solution variables $\lambda^n_m$ versus $n$. We can observe the new gaps that are created at each new generation: data are drawn as crosses, and lines join the $\lambda$ values of roots of $Z$ lying in the same gap for successive values of $n$. As a gap gets ``old'', the root inside it tends to a limit value. From Fig. \ref{porca1} we can also observe that all roots have absolute value less than $1/10$, while in principle they are bound by one.

\begin{figure}
\centerline{\includegraphics[width=.6\textwidth, angle = -90]{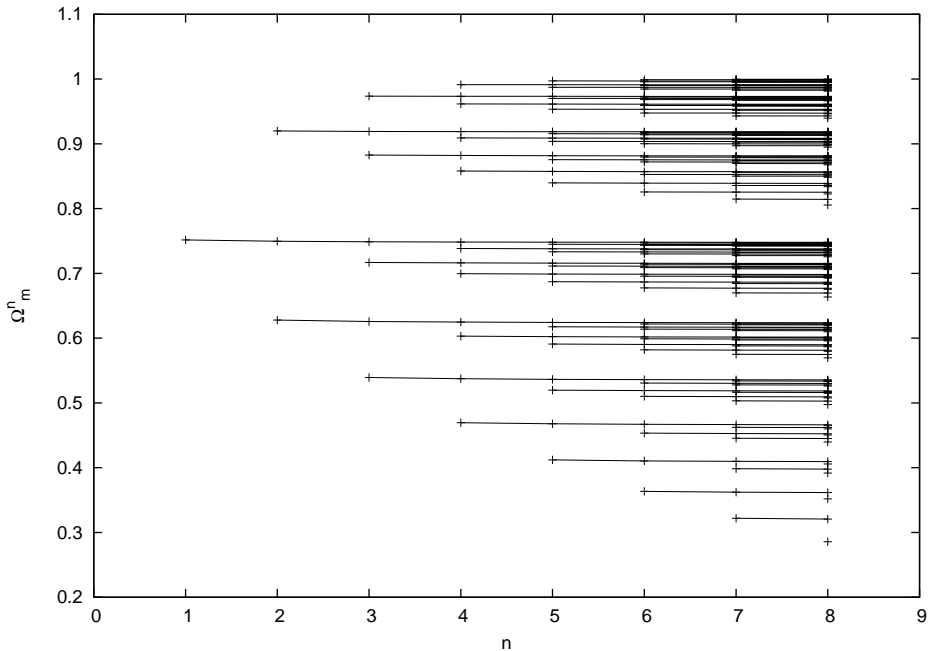}}
\caption{Integrated measures $\Omega^n_m$ versus versus generation number $n$, for the IFS described in the text. Lines connect points as explained in the text, and as done in the previous figure \ref{porca1}. They are only apparently horizontal (see the following Fig. \ref{porca4}). The number of Gaussian points employed is 160.}
\label{porca2}
\end{figure}

\begin{figure}
\centerline{\includegraphics[width=.6\textwidth, angle = -90]{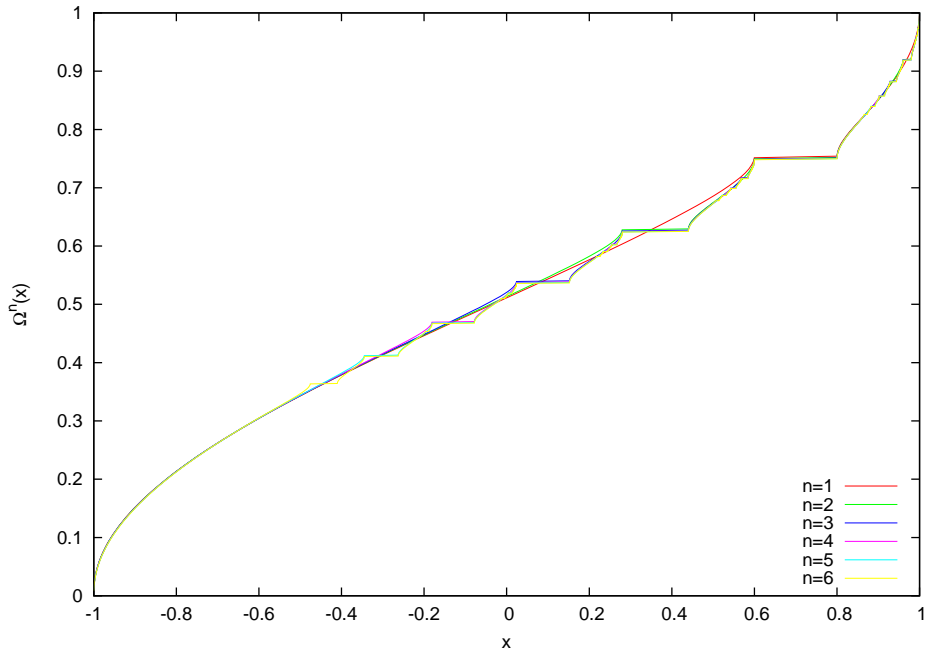}}
\caption{Harmonic measures $\Omega^n(x)$ versus $x$ at various generation numbers $n$  for the IFS described in the text. The number of Gaussian points employed is 160.}
\label{porca3}
\end{figure}

After all this preparatory work we are now ready to examine the equilibrium measure on the set $E^n$. It is straightforward to obtain the harmonic frequencies $\omega^{n}_i$ and the integrated measures $\Omega^{n}_mi$, discussed in the second section, via the Gaussian integration technique described in Sect. \ref{sec-integ}.
The case of fig. \ref{porca1} is examined again in figure \ref{porca2}, by plotting the integrated measures $\Omega^n_m$ versus $n$: recall that $\Omega^n_m$ is the measure of the set $(\gamma_1,\zeta^n_m)$ under the equilibrium measure on $E^n$. As such, these values are the heights of the plateaus in Fig. \ref{porca3}, that draws the integral
\begin{equation}
 \label{meas16}
   \Omega^n(x) = \int_{\gamma_1}^x d \sigma^n(s),
 \end{equation}
versus $x$, for the same IFS described above.

\begin{figure}
\centerline{\includegraphics[width=.6\textwidth, angle = -90]{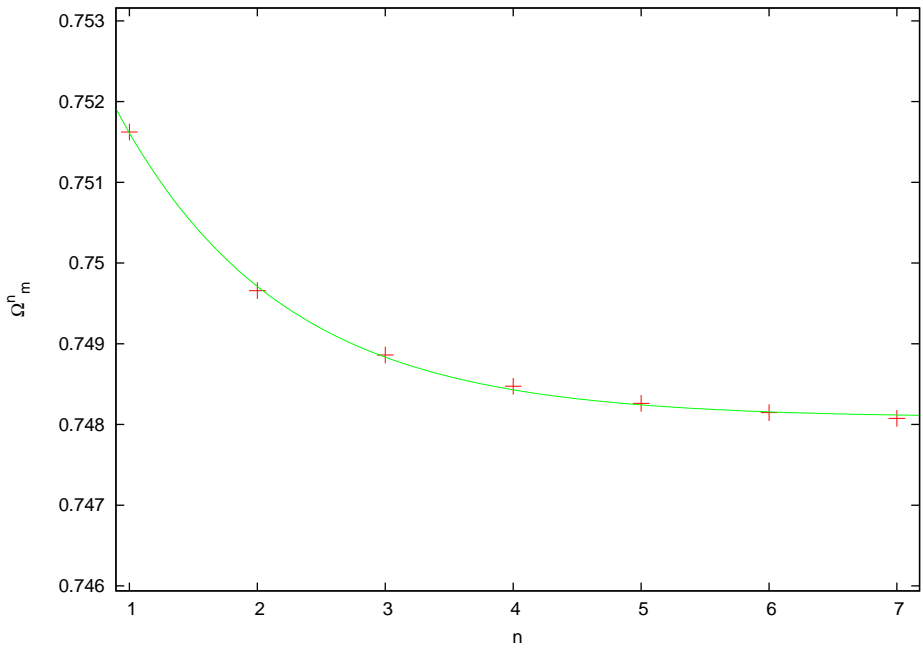}}
\caption{Integrated measure $\Omega^n_{2^n}$ versus generation number $n$ for the IFS described in the text (red crosses). Also plotted is the fit of the discrete values by the curve $f(n) = a + b e^{-cn}$ (green line). The number of Gaussian points employed is 160.}
\label{porca4}
\end{figure}

As $n$ grows, $E^n$ tends in Hausdorff distance to the attractor of the IFS, that is a Cantor set of zero Lebesgue measure. We want now to show that our technique can be used to compute the equilibrium measure of the IFS attractor itself: it appears from Fig. \ref{porca3} that $\Omega^n(x)$ tends to a limit function $\Omega(x)$ for Lebesgue a.e. $x$. This function is precisely the integral of the equilibrium measure on the IFS attractor.

In fact, we have already remarked that gaps, once created, stay forever in the complementary of the Cantor set. Let us therefore gauge the convergence, as $n$ grows, of the equilibrium measure of the infinite interval to the left of each gap. This is nothing else than $\Omega^n_{m(n)}$, when the index of the gap, $m$, depends on the generation $n$ according to the rule $m(n+1) = M m(n)$. In Figure \ref{porca4} we choose the first gap, $G^1_1$, that is labeled at successive generations as $G^n_{2^n}$: we extract a single line from fig. \ref{porca3} and we replot it alone so to render evident its variation. It appears that we have exponential convergence of the equilibrium measures, to the limit equilibrium measure of the Cantor set. We have found the same behavior for all gaps, so that we can numerically observe convergence to the limit function $\Omega(x)$ . But measures are born for integration, and integrating is what we are now up to.

\section{The Capacity of $E^n$ and of the IFS attractor}
\label{sec-capa}
The convergence that we have observed at the end of the previous section is just an instance of the convergence of $\sigma^n$ to the equilibrium measure on the IFS attractor. A second example is obtained considering the electrostatic potential $V(\sigma^n;z)$ defined in eq. (\ref{pote1}). Our technique permits to compute it in the full complex plane. Care must be exerted when $z$ belongs to $E^n$, for it may coincide with one of the Gaussian points, in which case one simply has to change their number in order to remove the coincidence. In this section we consider again the IFS generating the middle--third Cantor set, given by two maps with $\delta=1/3$ and $\gamma_1=-1$, $\gamma_2=1$.

\begin{figure}
\centerline{\includegraphics[width=.6\textwidth, angle = -90]{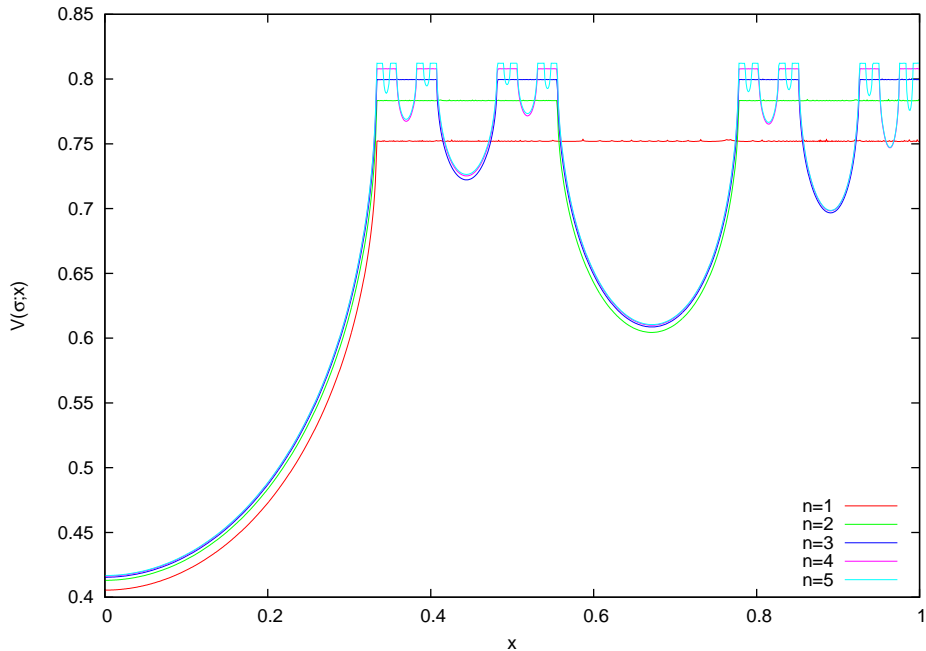}}
\caption{Electrostatic potential $V(\sigma_n;x)$ versus $x$, at generation level $n=1,\ldots,5$, for the IFS generating the ternary Cantor set. The number of Gaussian points employed is 2048.}
\label{fig-capa1}
\end{figure}

In Fig. \ref{fig-capa1} we plot the potential $V(\sigma^n;x)$ when $x$ belongs to the interval $[0,1]$, for $n=1,\ldots,5$. Notice that for symmetry we have constructed our ternary Cantor set so that its convex hull is $[-1,1]$. In the range of the figure, we compute the potential both on the set $E^n$ and on its complement. In the gaps, as the generation level $n$ grows, we observe exponential convergence to a limit value.
At fixed $n$, the potential $V(\sigma^n;x)$ must take a constant value almost everywhere on $E^n$. This is also observed in Fig. \ref{fig-capa1}. This constant value is linked to $C(E^n)$, the capacity of $E^n$, by
\begin{equation}
 \label{pote6}
    V(\sigma^n;x \in E^n)  = - \log( C(E^n)),\;\; \sigma^n \; \mbox{a.e.}
 \end{equation}
Our technique therefore also permits to compute these capacities, by computing the potential.

In Table 1 we plot $V(\sigma^n;x)$ versus $n$ for a value of $x$ in $E^n$, $n=7$, close to $\gamma_1$ (second column), together with the average potential, $\sum_l V(\sigma^n;x_l)/L$, taken over $L$ points in $E^7$ (third column). Observe that since $E^n \subset E^m$ for any $n \geq m$, this choice guarantees that the same set of sample points can be used for all values of $n$ in the Table. Difference between the two columns is due, and can serve to gauge, both the coarseness of Gaussian integration and the numerical errors involved in our technique. We observe a discrepancy that steadily diminishes as $n$ increases, from $2 \times 10^{-4}$ for $n=1$ to $3 \times 10^{-6}$ at $n=7$, in the same conditions.
\begin{table}
\centering
\begin{tabular}{|l|l|l|}
  \hline
  $n$ & $V(\sigma^n;x)$ & $\langle{V(\sigma^n;x)}\rangle$ \\
  1 & 0.751845  & 0.752051  \\
  2 & 0.799586 &  0.783380 \\
  3 & 0.807899 &  0.807941 \\
  4 & 0.807899 &  0.807941 \\
  5 & 0.812184 &  0.812210 \\
  6 & 0.814388 &  0.814392 \\
  7 & 0.815506 &  0.815509 \\
  \hline
\end{tabular}
\label{tab1}
\caption{Potential values $V(\sigma^n;x)$ at $x=-0.999996236647154$ and average values $\langle{V(\sigma^n;x)}\rangle$ over $L=4096$ points for the middle--third Cantor set IFS at different generation levels $n$. The number of Gaussian points employed is 2048.}
\end{table}

\begin{figure}
\centerline{\includegraphics[width=.6\textwidth, angle = -90]{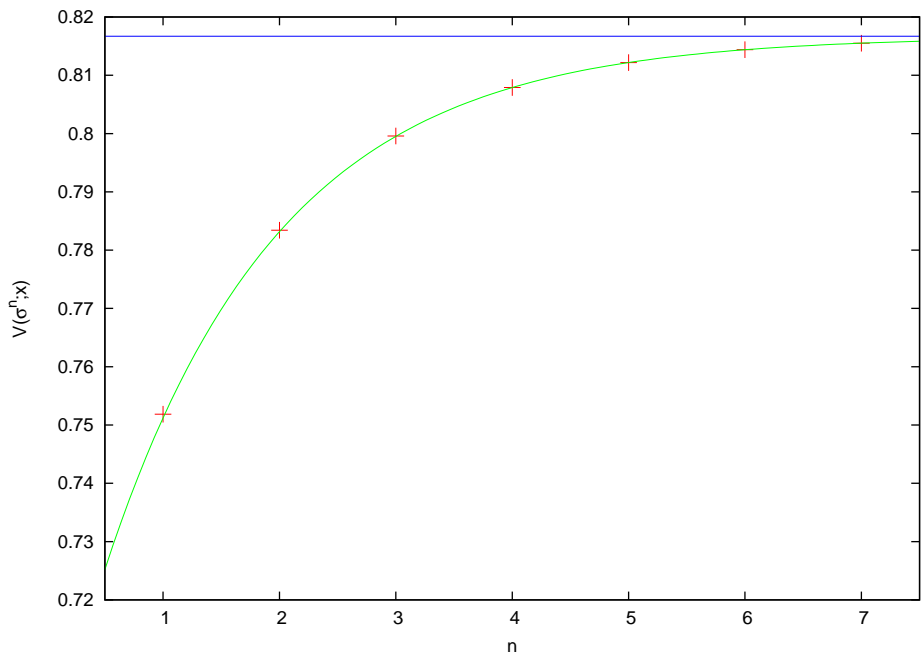}}
\caption{Data from the second column in Table I: potential values $V(\sigma^n;x)$ at $x=-0.999996236647154$ versus $n$ (red crosses). Also plotted are the fit of the discrete values from $n=4$ to $n=7$ by the curve $f(n) = a + b e^{-cn}$ (green line) and the infinite limit value $a=0.81668890$ (blue line). The other fitting parameters are $b=-0.1278376$ and $c=0.66927525$. }
\label{porca8}
\end{figure}

In Fig. \ref{fig-capa1} we observe convergence of the potential $V(\sigma^n;x)$ for increasing $n$ at fixed $x$, also for $x \in E^n$, although at a visibly slower rate than for $x$ in the complement of $E^n$.
This is particularly interesting if we now try to extrapolate the values in Table 1 to evaluate the capacity of the attractor of the IFS: Figure \ref{porca8} displays the same values of the second column of Table 1 together with a fit by an exponential function. Data from the third column would be indistinguishable from these, on the scale of the figure. Although much better extrapolation algorithms can be applied \cite{claudemichela}, this simple technique already provides an asymptotic value of the potential, which implies, via eq. (\ref{pote6}), an approximate value of the capacity $C({\cal A}) = 0.44189238$ that agrees within $6 \times 10^{-6}$ with the value $0.441898204379014$ provided by \cite{ran} (obviously doubled, because of the different scaling employed here). Performing the same extrapolation technique on the average values in the third column of Table 1 yields an even better estimate, $C({\cal A}) = 0.44189726$, this time within $1 \times 10^{-6}$ of the exact value. We can therefore conclude that we can reliably compute integrals with respect to the equilibrium measure on a Cantor set.


\section{Conclusions}

We have presented in this paper a technique for computing the equilibrium measure on attractors of iterated function systems. This technique consists of the careful concatenation of different ideas in numerical analysis, potential theory, and IFS construction. We have tested our approach by computing the capacity of the middle--third Cantor set, for which a reliable independent computation exists.

Yet, our interest in this problem does not end here: we intend to use the techniques of these paper as a mathematical electrophorus to compute the Jacobi matrix of the equilibrium measure. In fact, since this matrix embodies the orthogonal polynomials of the measure, its asymptotic properties are of paramount importance in a number of problems of importance in harmonic analysis, quantum mechanics and fractal geometry \cite{stahl,str0,str0b,poin1,poin2,etna,stric0,physd2}. We have reasonable hope that this investigation will help us to clarify a long-standing conjecture on the almost--periodicity of Jacobi matrices associated with IFS measures \cite{physd1}.

\end{document}